\documentclass{article}
\usepackage{amsfonts}
\usepackage[mathscr]{eucal}
 \usepackage{amsmath}
 \usepackage{amssymb}

\usepackage[T2A]{fontenc}
\usepackage[cp1251]{inputenc}
\usepackage[english]{babel}

\newenvironment{declaration}[1]{\trivlist
\item[\hskip \labelsep{\bf #1 }]\ignorespaces}{\endtrivlist}
\newenvironment{proofof}[1]{\begin{declaration}{#1}}{\hfill
$\square$ \end{declaration}}

\begin{document}

\title{Stability of a viral infection model with state-dependent delay, CTL and antibody immune responses}


\author{Alexander  Rezounenko
\\
 V.N.Karazin Kharkiv National University, \\ Kharkiv, 61022,  Ukraine\\
Institute of Information Theory and Automation,\\
Academy of Sciences  of the Czech Republic, \\ P.O. Box 18, 182\,08
Praha, CR\\E-mail:  rezounenko@yahoo.com }

  \maketitle
\begin{abstract}\marginpar{\tiny 2016.03.20.AR}

A virus dynamics model with intracellular state-dependent delay and
nonlinear infection rate of Beddington–DeAngelis functional response
is studied. The technique of Lyapunov functionals is used to analyze
stability of an interior infection equilibrium which describes the
case of both CTL and antibody immune responses activated. We
consider first a particular biologically motivated class of discrete
state-dependent delays. Next, the general case is investigated.

\par\noindent
{\bf Keywords: } evolution equations, Lyapunov stability,
state-dependent delay, virus infection model.
\par\noindent
{\bf 2010 MSC:} 
93C23, 
34K20, 
93D20 
97M60 

\end{abstract}

\section{Introduction}

We are interested in mathematical models of infectious diseases.
The 
diseases are caused by pathogenic microorganisms,
such as bacteria, viruses, parasites or fungi.
According to World Health Organization, many viruses (as  Ebola
virus, Zika virus, HIV, HBV, HCV and others) continue to be a major
global public health issues.


In our research we concentrate on models of viral infections.
 There are variety of models with and without delays which describe  dynamics between different viral infections and immune responses. Delays could be
concentrated or distributed. 
 We will not describe the historical evolution of such models, just
mention that early models
\cite{Perelson-Neumann-Markowitz-Leonard-Ho-S-1996,Nowak-Bangham-S-1996} contain three
variables: susceptible
host cells, infected cells and free virus. Next step was to take
into account, as written in \cite{Wodarz-2007-book}, that  "one
particular part of the immune system that is very important in the
fight against viral infections are the killer T cells or cytotoxic T
lymphocytes (CTL)." See also \cite{Zhu-Zou-DCDS-B-2009} and
references therein. There is another adaptive immune response by
antibodies.
The relative balance of both types of immune response "can be a decisive factor that determines whether patients are
asymptomatic or whether pathology is observe" \cite{Wodarz-JGV-2003}.
These lead to introduction of
 two additional variables of both immune responses
 \cite{Wodarz-JGV-2003,Wodarz-2007-book} (see also \cite{Yousfi-Hattaf-Tridane-JMB-2011} and references therein).

The model under consideration contains five variables: susceptible (noninfected) host cells $T$, infected cells $T^{*}$,
free virus $V$, a CTL response $Y$, and an antibody response $A$.
In case of bilinear nonlinearities and one constant concentrated
delay (see, for example, \cite{Yan-Wang-DCDS-B-2012}) it has the
following form
\begin{equation}\label{sdd-vir-01}
\left\{
\begin{array}{l}
   \dot T(t) = \lambda - d T(t) - k T(t)V(t), \\
 \dot T^{*}(t) = e^{-\omega h} k T(t-h)V (t-h) - \delta T^{*}(t) - p Y(t)T^{*}(t), \\
  \dot V(t) = N\delta T^{*}(t) - cV (t) - qA(t)V(t), \\
   \dot Y(t) = \beta T^{*}(t)Y(t) - \gamma Y(t)\\
   \dot A(t) = g A(t)V(t)  - bA(t). \\
\end{array}
\right.
\end{equation}
Here the dot over a function denotes the time derivative i.g, $\dot
T(t) = {dT(t)\over dt}$, all the constants $\lambda, d, k, \delta,
p, N, c, q, \beta, \gamma, g, b,\omega$ are positive.
As for the immune responses, the fourth equation describes the regulation of CTL response with $p Y(t)T^{*}(t)$ (in the second equation) being the rate of killing of infected cells by lytic immune response. The fifth equation describes the regulation of antibody  response with $qA(t)V(t)$ (in the third equation) being the rate of virus neutralization  by antibodies \cite[p.1744]{Wodarz-JGV-2003}.
 In (\ref{sdd-vir-01}), $h$ denotes the delay between the time the virus contacts a target cell and the
time the cell becomes actively infected (starts to produce new virions).

In the above model, the standard bilinear incidence rate is used according to the principle of mass action.
For more details and references on the models of infectious diseases
with more general types of nonlinear incidence rates see e.g. \cite{Korobeinikov-BMB-2007,Gourley-Kuang-Nagy-JBD-2008}.
In  paper \cite{Zhao-Xu_EJDE-2014}, following
\cite{Huang-Ma-Takeuchi-AML-2011,Wang-Liu-MMAS-2013,Yan-Wang-DCDS-B-2012},  authors assume that
the infection rate of the virus dynamics models is given by the
Beddington-DeAngelis functional response
\cite{Beddington-JAE-1975,DeAngelis-Goldstein-ONeill-E-1975},
${kTV\over 1+k_1 T+k_2 V}$, where $k,k_1,k_2~\ge~0$ are constants.
The Lyapunov asymptotic stability\cite{Lyapunov-1892} of points of
equilibrium is studied for the following model with {\it constant}
concentrated delay
\begin{equation}\label{sdd-vir-02}
\left\{
\begin{array}{l}
   \dot T(t) = \lambda - d T(t) - f(T(t),V(t)), \\
 \dot T^{*}(t) = e^{-\omega h} f(T(t-h),V (t-h)) - \delta T^{*}(t) - p Y(t)T^{*}(t), \\
  \dot V(t) = N\delta T^{*}(t) - cV (t) - qA(t)V(t), \\
   \dot Y(t) = \beta T^{*}(t)Y(t) - \gamma Y(t)\\
   \dot A(t) = g A(t)V(t)  - bA(t). \\
\end{array}
\right.
\end{equation}
The functional response
\begin{equation}\label{sdd-vir-03}
f(T,V)={kTV\over 1+k_1 T+k_2 V}, \qquad k, k_1, k_2 \ge 0, \quad
T,V\in R,
\end{equation}
was introduced by Beddington \cite{Beddington-JAE-1975} and
DeAngelis et al.~\cite{DeAngelis-Goldstein-ONeill-E-1975}.

It is evident that the constancy of the delay is an extra assumption
which essentially simplifies the analysis, but is not motivated by
the biological background of the model. It was a reason (see e.g. \cite{Wang-Pang-Kuniya-Enatsu_AMC-2014})
to discuss distributed delay models as an alternative to discrete constant delay ones. We propose another approach.

Our goal is to remove the restriction of the constancy of the delay
and investigate the well-posedness and Lyapunov stability of the
following virus infection model with Beddington-DeAngelis
 functional response and {\it state dependent delay}.
It appears that the analysis essentially differs  from the constant delay case.
 To the best of our knowledge, such models have not been considered before.
It is well known that differential equations with state dependent delay are always nonlinear by its nature (see the review \cite{Hartung-Krisztin-Walther-Wu-2006} for more details and discussion).

  As usual in a delay system with
 (maximal) delay $h>0$ \cite{Hale,Kuang-1993_book,Walther_book}, for a function $v(t), t\in [a,b]\subset R, a>b+h$,
 we denote the history segment $v_t=v_t(\theta)\equiv v(t+\theta), \theta\in [-h,0].$ We
 denote the space of continuous functions by $C\equiv C([-h,0]; R^5)$ equipped with the sup-norm.
  In the above notations, we use $u(t)=(T(t),T^{*}(t),V(t),Y(t),A(t))$ and consider a continuous
 functional (state dependent delay) $\eta : C \to [0,h]$. Now we can
 present the system under consideration
\begin{equation}\label{sdd-vir-04}
\left\{
\begin{array}{l}
   \dot T(t) = \lambda - d T(t) - f(T(t),V(t)), \\
 \dot T^{*}(t) = e^{-\omega h} f(T(t-\eta(u_t)),V (t-\eta(u_t))) - \delta T^{*}(t) - p Y(t)T^{*}(t), \\
  \dot V(t) = N\delta T^{*}(t) - cV (t) - qA(t)V(t), \\
   \dot Y(t) = \beta T^{*}(t)Y(t) - \gamma Y(t)\\
   \dot A(t) = g A(t)V(t)  - bA(t). \\
\end{array}
\right.
\end{equation}
with the functional response $f(T,V)$ given by (\ref{sdd-vir-03}).
We denote by ${\cal F} $ the right-hand side of (\ref{sdd-vir-04})
to write the system shortly as $\dot u(t) = {\cal F} (u_t)$.

The paper is organized as follows. In Section 2, we discuss and
choose different sets of initial data and prove the existence and uniqueness of solutions.
 Next we prove that the sets are invariant. Section 3 is devoted to the stability properties of a
 stationary solution. We study the stability of an interior equilibrium which describes the case when both
 CTL and antibody immune responses are activated. 
We believe this infection equilibrium is only biologically meaningful in the study of the disease.
 We use the technique of Lyapunov functionals \cite{Lyapunov-1892} and consider first
 a particular (biologically motivated) case of state-dependent delay.
 We  complete the paper by investigation of the general case.

\section{Preliminaries}


We first study the basic questions of the existence and uniqueness
of solutions to the problem (\ref{sdd-vir-04}). Since two functions
$T$ and $V$ are used in (\ref{sdd-vir-04}) at different time moments
(current time $t$ and delayed time $t-\eta(u_t)$), we should
consider initial values $T(\theta), V(\theta)$ for $\theta\in
[-h,0]$. As usual for such a biological system, one should check the
non-negativeness and boundedness of all the coordinates provided
initial values are non-negative.

We will study the system  (\ref{sdd-vir-04}) with an initial
function


\begin{equation}\label{sdd-vir-19}
u_0 = \varphi \equiv (T_0,T^{*}_0,V_0,A_0,Y_0) \in C_{+}\equiv C_{+}[-h,0],
\end{equation}
where $R_{+}\equiv [0,+\infty), C_{+}\equiv C_{+}[-h,0]\equiv C([-h,0];
R^5_{+})$.

Let us introduce the  set
$$\Omega_C \equiv \left\{ \varphi \equiv (T_0,T^{*}_0,V_0,A_0,Y_0) \in C_{+}[-h,0], \right.
$$
$$ 0\le T_0(\theta)\le {\lambda \over d},\quad  0\le T^{*}_0(\theta)\le {k \lambda \over d k_2
\delta}e^{-\omega h}, \quad 0\le V_0(\theta)\le {N k \lambda \over c  d
k_2} e^{-\omega h},
$$
\begin{equation}\label{sdd-vir-20}
\left. 0\le T^{*}_0(\theta)+ {p\over \beta} Y_0(\theta) \le {k^2
\lambda^2 e^{-2 \omega h} \over d^2 c k_2 \min \{ \delta; \gamma\}}, \quad
0\le V_0(\theta)+ {q\over g} A_0 (\theta) \le {N k \lambda \, e^{-\omega h}
\over d k_2 \min \{ c; b\}}, \, \theta \in [-h,0] \right\}.
\end{equation}

\smallskip

We consider the following assumption on the state-dependent delay 
\begin{equation}\label{sdd-vir-H1}
{\bf (H1_\eta)}\quad\forall \psi \in Z^{2,3}\equiv \left\{
\psi=(\psi^1,\psi^2,\psi^3,\psi^4,\psi^5) \in C_{+} :
\psi^2(0)=\psi^3(0)=0\right\} \quad \Longrightarrow 
\eta(\psi)>0.
\end{equation}

\medskip

{\bf Remark 1.} {\it We notice that even more restrictive assumption
$\eta (\psi)>0$ for all $\psi\in C_{+}$ is biologically well motivated.
On the other hand, even this restriction (the so-called "non-vanishing delay") does {\tt not}
guarantee the uniqueness of solutions with merely continuous data (see \cite{Driver-AP-1963}).
}

\smallskip

The first result is the following

\medskip

{\bf Theorem 2.} {\it Let $\eta : C \to [0,h]$ be a continuous
 functional (state dependent delay). Then

(i) for any initial function $\varphi \in C
$ there exist {\tt continuous} solutions 
 to (\ref{sdd-vir-04}), (\ref{sdd-vir-19}).

 (ii) If additionally, $\eta$ satisfies $(H1_\eta)$, then for any initial function $\varphi \equiv (T_0,T^{*}_0,V_0,A_0,Y_0) \in \Omega_C
$ such that $T_0, V_0$ are Lipschitz functions, the problem (\ref{sdd-vir-04}),
(\ref{sdd-vir-19}) has a unique solution. The solution is globally
Lipschitz in time and satisfies
$$ u_t\equiv (T_t,T^{*}_t,V_t,A_t,Y_t)\in \Omega_C, \quad t \ge 0.
$$
}

%
%

\smallskip

{\it Proof of Theorem 2.} (i) The existence of continuous solutions
is guaranteed by the continuity of the right-hand side of
(\ref{sdd-vir-04}) and classical results on delay equations
\cite{Hale,Walther_book}.

(ii) Since $T_0, V_0$ are Lipschitz functions,  the uniqueness of
continuous solutions follows from the general results on differential equations with
state-dependent delay (see the review on ordinary equations
\cite{Hartung-Krisztin-Walther-Wu-2006} for details and references and also \cite{Rezounenko_NA-2009,Rezounenko_NA-2010,Rezounenko_JMAA-2012,Rezounenko-Zagalak-DCDS-2013} for PDEs)).
Let us show that the set $\Omega_C$ is invariant i.e. any solution
starting from $\varphi\in \Omega_C$ remains in $\Omega_C$ for all
$t\ge 0.$

We notice that in the case of {\it constant} delay  the non-negativeness of all coordinates of a solution
follows from the {\it
quasi-positivity} property of the right-hand side of
(\ref{sdd-vir-04}) (see e.g. \cite[Theorem 2.1, p.81]{Smith-1995-book}).
We stress that in the case of {\it state-dependent} delay we cannot directly apply
\cite[Theorem 2.1, p.81]{Smith-1995-book} because it relies on the
Lipschitz property of the right-hand side of a system, which is {\it
not} the case for (\ref{sdd-vir-04}). We could use the corresponding
extension to the state-dependent delay case
\cite{Rezounenko-JADEA-2012}, but we propose another way here.

To prove the non-negativeness of all coordinates  of a solution
$u(t)=(T(t),T^{*}(t),V(t),Y(t),A(t))$ we use the direct analysis of each coordinate.
It is easy to see that $T(t)\to 0+$ implies $\dot T(t) \to \lambda >0$ which makes
impossible for $T$ to become negative.  The direct integration shows that coordinates satisfy
\begin{equation}\label{sdd-vir-31}
T^{*}(t)= T^{*}(0) e^{-\int^t_0 (\delta+pY(s))\, ds}
+ e^{-\omega h}\int^t_0 f(T(\tau-\eta(u_\tau)),V (\tau-\eta(u_\tau)))
 e^{-\int^t_\tau (\delta+pY(s))\, ds} \, d\tau,
\end{equation}
\begin{equation}\label{sdd-vir-32}
V(t)=V(0)\, e^{-\int^t_0 (c+qA(s))\, ds}
+ N\delta \int^t_0 T^{*}(\tau)
 e^{-\int^t_\tau (c+qA(s))\, ds} \, d\tau,
\end{equation}
\begin{equation}\label{sdd-vir-33}
Y(t)=Y(0)\, e^{-\int^t_0 (\beta T^{*}(\tau)-\gamma))\, d\tau}, \qquad
A(t)=A(0)\, e^{-\int^t_0 (g V(\tau)-b))\, d\tau}.
\end{equation}
Equations (\ref{sdd-vir-33}) show that $Y(0)\ge 0, A(0)\ge 0$ implies $Y(t)\ge 0, A(t)\ge 0$ for all $t\ge 0$.
For the constant delay case, equations (\ref{sdd-vir-31}),  (\ref{sdd-vir-32}) would imply the similar
result for $T^{*}(t), V(t)$, but for the {\it state-dependent} delay we need more care. First, (\ref{sdd-vir-31}) shows the property (for some $t^1\ge 0$)
\begin{equation}\label{sdd-vir-34}
V(s)\ge 0,\, s\in [-h,t^1] \quad \mbox{ implies } \quad T^{*}(s)\ge 0, \, s\in [0,t^1].
\end{equation}
Similarly, (\ref{sdd-vir-32}) gives
\begin{equation}\label{sdd-vir-35}
T^{*}(s)\ge 0, \, s\in [0,t^1] \quad \mbox{ implies } \quad V(s)\ge 0, \, s\in [0,t^1].
\end{equation}
Now, let us assume that the non-negativity of $T^{*}$ or $V$ falls. Properties (\ref{sdd-vir-34}), (\ref{sdd-vir-35}) show
that $T^{*}$ and $V$ should change the sign simultaneously i.e. there exist a (smallest possible) time moment $t^1\ge 0$
and $\delta^1>0$ such that $T^{*}(t)\ge 0, V(t)\ge 0$ for $t\in (t^1-\delta^1, t^1]$ and
$T^{*}(t)< 0, V(t)< 0$ for $t\in (t^1, t^1+\delta^1)$. By the continuity of solutions, $T^{*}(t^1)=V(t^1)=0$ which implies
(see $(H1_\eta)$) that $u_{t^1}\in Z^{2,3}$ i.e. $\eta(u_{t^1})>0$. Hence there exists $\delta^2>0$ such that $V (\tau-\eta(u_\tau))\ge 0$ for
$\tau\in (t^1, t^1+\delta^2)$. By this property and (\ref{sdd-vir-31}), one has  $T^{*}(t)\ge 0$ for  $t\in (t^1, t^1+\delta^2)$.
It contradicts the choice of $t^1$ and completes the proof of the non-negativity of all coordinates.

Let us prove the upper bounds on the
coordinates, given in (\ref{sdd-vir-20}). To save the space we
formulate an easy variant of the Gronwall's lemma.

 \medskip

 {\bf Lemma 3.} {\it Let $\ell \in C^1 [a,b)$ and satisfy  ${d\over dt}
\ell (t) \le c_1 - c_2 \ell (t), \, t\in [a,b)$. Then $\ell (a)\le
c_1 c^{-1}_2 $ implies $\ell (t)\le c_1 c^{-1}_2 $ for all $t\in
[a,b)$. In the case $b=+\infty$, for any $\varepsilon>0$ there
exists $t_\varepsilon \ge a$ such that $\ell (t)\le c_1 c^{-1}_2 +
\varepsilon$ for all $t\ge t_\varepsilon$.
}

{\it Proof of Lemma 3} is easy multiplication of the inequality
by $e^{c_2 t}$ and integration over $[a,t]$. It leads to $\ell (t)
\le \left( \ell (a) - {c_1\over c_2}\right) e^{-c_2(t-a) }+
{c_1\over c_2}.$ It completes the proof of Lemma 3. $\Box $

\smallskip

Since $f$ is non-negative  for non-negative arguments (see
(\ref{sdd-vir-03})), we get from the first equation of
(\ref{sdd-vir-04}) the estimate $\dot T(t) \le \lambda - d T(t)$.
Hence Lemma 3 and $T(0)\le {\lambda \over d}$ implies $T(t)\le
{\lambda \over d}$ for $t \ge 0$. We use it to estimate the second
coordinate, see (\ref{sdd-vir-03}), as follows $ f(T,V) \le {k
\lambda V\over d (1+ k_2 V)} \le {k \lambda \over d k_2}$. It gives
$\dot T^{*}(t) \le {k \lambda \over d k_2}e^{-\omega h} - \delta T^{*}(t)$
and Lemma 3 implies the needed bound for $T^{*}$ in
(\ref{sdd-vir-20}). The bound for $T^{*}$ and the third equation in
(\ref{sdd-vir-04}) give $\dot V(t) \le N\delta T^{*}(t) - c V(t) \le
{N  k \lambda \over d k_2}e^{-\omega h}  - c V(t)$. Lemma 3 proves the
estimate for $V$ in (\ref{sdd-vir-20}). Next, we use the second and
the fourth equations in (\ref{sdd-vir-04}) to get
$$ \dot T^{*}(t) + {p\over \beta} \dot Y(t)= e^{-\omega h} f(T(t-\eta(u_t)),V (t-\eta(u_t))) - \delta T^{*}(t)  - {\gamma p\over \beta} Y(t)
$$
$$\le e^{-\omega h} k T_{max} V_{max} - \min \{ \delta; \gamma \} \left( T^{*}(t) + {p\over \beta} Y(t)\right)
$$
$$\le  k {\lambda \over d}
{N k \lambda \over c  d
k_2} e^{-2\omega h} 
- \min \{ \delta; \gamma \} \left( T^{*}(t) + {p\over \beta}
Y(t)\right).
$$
Lemma 3 proves the bound for $T^{*}(t) + {p\over \beta} Y(t)$ in
(\ref{sdd-vir-20}). In the similar way, using the third and fifth
equations in (\ref{sdd-vir-04}), one gets
$$ \dot V(t)+{q\over g} \dot A(t) \le N\delta T^{*}(t) - cV (t) - {bq\over g} A(t)
    \le {N k \lambda \over d k_2}e^{-\omega h} - \min \{ c; b\} \left( V(t)+ {q\over g} A(t)
    \right).
$$
Lemma 3 implies the last estimate in (\ref{sdd-vir-20}). All solutions are global (defined for all $t\ge -h$).
It completes the proof of Theorem 2. $\Box $

\smallskip

{\bf Remark 4.} {\it We notice the our invariant set $\Omega_C$
differs from the absorbing set $\Gamma$, used in
\cite{Zhao-Xu_EJDE-2014} for the constant delay system. Let us
denote by $\Omega^\varepsilon_C$ the set where all the upper bounds
in (\ref{sdd-vir-20}) are increased by $\varepsilon$. Then the
second part of the Lemma 3 implies that for any $\varepsilon>0$ the
set $\Omega^\varepsilon_C$
is absorbing for any solution (not necessary starting in
$\Omega_C$). Another difference is that all the five coordinates of
$\varphi\in \Omega_C$ are continuous functions in contrast to the
constant delay case \cite{Zhao-Xu_EJDE-2014}, where the second,
fourth and fifth coordinates belong to $R_{+}$.}

\smallskip

{\bf Remark 5.} {\it It is well known that continuous solutions to
differential equations with state-dependent delay may be non-unique
(see examples in \cite{Driver-AP-1963}). There are two ways to
insure the uniqueness of solutions. The first one is to restrict the
set of initial functions \cite{Hartung-Krisztin-Walther-Wu-2006}.
The second one is to restrict the class of state-dependent delays
\cite{Rezounenko_NA-2009,Rezounenko_JMAA-2012} and work with
continuous initial functions.}

\smallskip

If one is interested in  continuously differentiable solutions, we
could also apply the solution manifold approach
\cite{Walther_JDE-2003,Hartung-Krisztin-Walther-Wu-2006} to the
initial value problem (\ref{sdd-vir-04}), (\ref{sdd-vir-19}). We
remind the short notation for the system (\ref{sdd-vir-04}) as $\dot
u(t) = {\cal F} (u_t)$.
 Let us introduce the following subset of $\Omega_C$ (c.f.
(\ref{sdd-vir-20}))
\begin{equation}\label{sdd-vir-21}
\Omega_{\cal F} \equiv \left\{ \varphi = (T_0,T^{*}_0,V_0,A_0,Y_0)
\in C^1_{+}[-h,0];\quad \varphi \in \Omega_C;\quad  \dot \varphi(0)
= {\cal F} (\varphi) \right\}.
\end{equation}

The following result is a corollary of Theorem 2.

\medskip

{\bf Theorem 6.} {\it Let $\eta : C \to [0,h]$ be a continuous
 functional (state dependent delay), 
 satisfying $(H1_\eta)$. Then
for any initial function $\varphi \in \Omega_{\cal F}
$ there exists a {\tt unique} ({continuously differentiable}) solution 
 to (\ref{sdd-vir-04}),
(\ref{sdd-vir-19}), satisfying
$$ u_t\equiv (T_t,T^{*}_t,V_t,A_t,Y_t)\in \Omega_{\cal F}, \quad t \ge 0.
$$
}

\smallskip

One can see that any continuous solution $u$ started at $\varphi \in
\Omega_C$, satisfies $u_t \in \Omega_{\cal F} $ for $t>h$.

\subsection{Stationary solutions}

We look for stationary solutions of (\ref{sdd-vir-04}).  Consider
the system with $u(t)=u(t-\eta(u_t))=\widehat u$ and denote the
coordinates of a stationary solution by
$(\widehat{T},\widehat{T^{*}},\widehat{V},\widehat{Y},\widehat{A})=\widehat
u\equiv \widehat{\varphi}(\theta),\, \theta\in [-h,0]$. We have
\begin{equation}\label{sdd-vir-stationary-1}
 \left\{
\begin{array}{l}
   0 = \lambda - d \widehat{T} - f(\widehat{T},\widehat{V}), \\
 0 = e^{-\omega h} f(\widehat{T},\widehat{V}) - \delta \widehat{T^{*}} - p \widehat{Y}\widehat{T^{*}}, \\
  0 = N\delta \widehat{T^{*}} - c \widehat{V} - q \widehat{A}\widehat{V}, \\
   0 = \beta \widehat{T^{*}}\widehat{Y} - \gamma \widehat{Y}\\
   0 = g \widehat{A}\widehat{V}  - b\widehat{A}. \\
\end{array}
\right.
\end{equation}

Since the stationary solutions of
(\ref{sdd-vir-04}) are the same as
for a constant delay system, we proceed as follows 
(see e.g. \cite{Zhao-Xu_EJDE-2014}). The last two equations in
(\ref{sdd-vir-stationary-1}) give $\widehat{T^{*}}= {\gamma \over
\beta }, \widehat{V} = {b\over g}$. Using this, the third equation
implies $\widehat{A} = { N\delta\gamma g - \beta cb\over \beta
qb}.$  To insure the positivity of $\widehat{A}$ we assume the constants in the system (\ref{sdd-vir-04}) satisfy

\begin{equation}\label{sdd-vir-H2}
\hskip-70mm {\bf (H2)} \hskip50mm N\delta\gamma g > \beta cb.
\end{equation}

Now we substitute the value
$\widehat{V} = {b\over g}$ into the first equation of
(\ref{sdd-vir-stationary-1}) and use (\ref{sdd-vir-03}) to get the
quadratic equation for $\widehat{T}$
\begin{equation}\label{sdd-vir-stationary-3} 
 dgk_1 \widehat{T}^2   + (dk_2 b + dg - \lambda gk_1 + k b)\, \widehat{T} - \lambda (g+k_2b) = 0.
\end{equation}
One can easily see that the discriminant of the above equation is
positive and there are two roots of different signs. We denote below by
$\widehat{T}$ the unique {\it positive} root.  The first two
equations in (\ref{sdd-vir-stationary-1}) give (remind that
$\widehat{T^{*}}$ is already known) $\widehat{Y} = {\lambda - d
\widehat{T} - e^{\omega h} \delta \widehat{T^{*}}\over e^{\omega h} p
\widehat{T^{*}}}$. The positivity of $\widehat{Y}$ follows from
the following assumption

\begin{equation}\label{sdd-vir-H3}
\hskip-50mm {\bf (H3)} \hskip50mm \quad
\lambda > d \widehat{T} +
{\delta \gamma
\beta^{-1}} e^{\omega h},
\end{equation}
where $\widehat{T}$ is the unique {\it positive} root
of the quadratic equation (\ref{sdd-vir-stationary-3}).

\smallskip

We notice that, from biological point of view, (H2), (H3) are standard assumptions on reproduction numbers, which
are given here in a short form.
We could summarize the above calculations in the following

\medskip

{\bf Lemma 7.} 
{\it Let assumptions (H2) and (H3) be satisfied (see (\ref{sdd-vir-H2}), (\ref{sdd-vir-H3})).
Then the system (\ref{sdd-vir-stationary-1})
has a unique solution $(\widehat{T},\widehat{T^{*}},\widehat{V},\widehat{Y},\widehat{A})$ (the unique stationary solution of
(\ref{sdd-vir-04})). All the coordinates are positive, $\widehat{T}$ is the unique {\it positive} root
of the quadratic equation (\ref{sdd-vir-stationary-3}) and coordinates satisfy
\begin{equation}\label{sdd-vir-stationary-2}
\left\{
\begin{array}{l}
 \widehat{T^{*}}= {\gamma \over \beta }, \quad \widehat{V} = {b\over
   g}, \quad    \widehat{A} = { N\delta\gamma g - \beta cb\over \beta qb},  \quad \widehat{Y} = {\lambda - d
\widehat{T} - e^{\omega h} \delta \widehat{T^{*}}\over e^{\omega h} p
\widehat{T^{*}}},\\
  N\delta \widehat{T^{*}} =  \widehat{V} (c+ q \widehat{A}), \quad
  \lambda = d \widehat{T} + f(\widehat{T},\widehat{V}),\quad  (\delta  - p \widehat{Y})\widehat{T^{*}}e^{\omega h}=  f(\widehat{T},\widehat{V}). \\
\end{array}
\right.
\end{equation}
}

\smallskip

These equations connecting the coordinates of the stationary
solution will be used below in the study of the stability properties
(c.f. (\ref{sdd-vir-stationary-1})).

\section{Stability properties}

The function $v(x) = x - 1- \ln x$ for $x>0$ plays an important role
in 
construction of Lyapunov functionals. One can
easily check that $v(x)\ge 0$ and $v(x)=0$ if and only if $x=1$. The
derivative equals $\dot v(x) = 1-{1\over x}$, which is evidently
negative for $x\in (0,1)$ and positive for $x>1$. The graph of $v$
explains the use of the composition $v\left({x\over x^0}\right)$ in
the study of the stability properties of an equilibrium $x^0$.
Another important property is the following estimate
\begin{equation}\label{sdd-vir-v-functional-2}
\forall \delta \in (0,1)\quad \forall x \in (1-\delta, 1+\delta)
\quad \mbox{ one has } \quad {(x-1)^2\over 2(1-\delta)} \le v(x) \le
{(x-1)^2\over 2(1+\delta)}.
\end{equation}
To check it, one simply observes that all three functions vanish at
$x=1$ and $\left|{d\over dx} \left( {(x-1)^2\over 2(1-\delta)}
\right) \right| \le  |{d\over dx}  v(x)| \le \left|{d\over dx}
\left( {(x-1)^2\over 2(1+\delta)} \right) \right|$ in the
$\delta$-neighborhood of $x=1$.

\subsection{Particular case of a state-dependent delay}

As before, we denote $u(t)=(T(t),T^{*}(t),V(t),Y(t),A(t))$. Consider
arbitrary $\varphi =
(\varphi^1,\varphi^2,\varphi^3,\varphi^4,\varphi^5)\in C$. We are
interested in the following particular form of the state-dependent
delay
\begin{equation}\label{sdd-vir-05}
\eta (\varphi)=F(\varphi^1(0),\varphi^3(0)).
\end{equation}
It means that
 $\eta (u_t) = F(T(t),V(t))$ which looks natural since the delay appears in the nonlinearity
 $f$ which depends on $T$ and $V$ only (see the second equation in (\ref{sdd-vir-04})).

\medskip

{\bf Theorem 8}. {\it Let assumptions (H2) and (H3) be satisfied (see (\ref{sdd-vir-H2}), (\ref{sdd-vir-H3})).
Assume the state-dependent delay $\eta$ has the form
(\ref{sdd-vir-05}) with a continuous $F: R^2_{+}\to [0,h]$,
satisfies $(H1_\eta)$ (see (\ref{sdd-vir-H1})) and
\begin{equation}\label{sdd-vir-06}
|\eta (\varphi)-\eta (\widehat\varphi)|
 = |F(\varphi^1(0),\varphi^3(0))-F(\widehat{T},\widehat{V})|
 \le c_\eta \left( (\varphi^1(0)-\widehat{T})^2 + (\varphi^3(0)-\widehat{V})^2\right).
\end{equation}
Then the stationary solution
$\widehat{\varphi}=(\widehat{T},\widehat{T^{*}},\widehat{V},\widehat{Y},\widehat{A})$
is locally asymptotically stable. For sufficiently small values of $c_\eta$, the stationary
solution is globally asymptotically  stable. }

\smallskip

 {\it Proof of Theorem 8.}
Consider the 
Lyapunov functional
$$U^{1}(t) \equiv  \left( T(t) - \widehat T -
\int^{T(t)}_{\widehat T} {f(\widehat{T},\widehat{V})\over
f(\theta,\widehat{V})} \, d\theta  \right)e^{-\omega h}
 + \widehat{T^{*}}\cdot v\left({T^{*}(t)\over \widehat{T^{*}}}\right) +
{\delta + p\widehat{Y}\over N\delta} \cdot v\left({V(t)\over
\widehat V}\right) $$
\begin{equation}\label{sdd-vir-Lyapuniv-functional-1}
+ {p\over \beta} \cdot v\left({Y(t)\over \widehat Y}\right) + 
{q\over Ng}\left( 1+ {p\widehat{Y}\over \delta}\right) \cdot
v\left({A(t)\over \widehat A}\right) + (\delta + p\widehat{Y})
\widehat{T^{*}} \int^t_{t-\eta(\widehat{\varphi})}
v\left({f(T(\theta),V(\theta))\over f(\widehat{T},\widehat{V})}
\right) \, d\theta.
\end{equation}
 We use the same notations as in \cite{Zhao-Xu_EJDE-2014} to simplify
for the reader the comparison of the computations. In spite of the
same Lyapunov functional as in \cite{Zhao-Xu_EJDE-2014}, the time
derivative of $U^{1}(t)$ along a solution $u$ of (\ref{sdd-vir-04})
is different due to the state-dependence of the delay in the system.
It reads as follows
$${d\over dt}U^{1}(t) = \left( 1- {f(\widehat{T},\widehat{V})\over f({T}(t),\widehat{V})}\right)
 e^{-\omega h} \left( \lambda - d T(t) - f(T(t),V(t))\right)
$$
$$ + \left( 1 - { \widehat{T^{*}}\over T^{*}(t)}\right)\left( e^{-\omega h} f(T(t-\eta(u_t)),V (t-\eta(u_t))) - \delta T^{*}(t) - p Y(t)T^{*}(t)\right)
$$
$$+ {\delta+P \widehat{Y}\over N\delta} \left( 1 - { \widehat{V}\over V(t)}\right)\left( N\delta T^{*}(t) - cV (t) - qA(t)V(t)\right)
$$
$$+ {p\over \beta} \left( 1 - { \widehat{Y}\over Y(t)}\right)\left( \beta T^{*}(t)Y(t) - \gamma Y(t)\right) + {q\over Ng}
\left( 1+ {p \widehat{Y}\over \delta}\right)\left( 1 - { \widehat{A}\over A(t)}\right)\left( g A(t)V(t)  - bA(t)\right)
$$
$$ + e^{-\omega h} \left[f(T(t),V (t)) - f(T(t-\eta(\widehat{\varphi})),V
(t-\eta(\widehat{\varphi}))) 
\right] + \widehat{T^{*}}(\delta + p\widehat{Y})
 \ln { f(T(t-\eta(\widehat{\varphi})),V
(t-\eta(\widehat{\varphi}))) \over f(T(t),V (t))}.
$$
Opening parenthesis, grouping similar terms and canceling some of
them, we obtain
$${d\over dt}U^{1}(t) =  \left( 1- {f(\widehat{T},\widehat{V})\over f({T}(t),\widehat{V})}\right)
 e^{-\omega h} d \left( \widehat{T} -  T(t)\right)
$$
$$ - \widehat{T^{*}}(\delta + p\widehat{Y}) \left[ {f(\widehat{T},\widehat{V})\over f({T}(t),\widehat{V})} -
{f({T}(t),{V}(t))\over f({T}(t),\widehat{V})} + {e^{-\omega h}\over \delta
+ p\widehat{Y}} \cdot  { f(T(t-\eta(u_t)),V (t-\eta(u_t))) \over
T^{*}(t)} \right. $$
$$ + \left. { T^{*}(t)\cdot \widehat{V}\over \widehat{T^{*}}\cdot  V(t)} +
{ V(t)\over \widehat{V} } - 3 - \ln {
f(T(t-\eta(\widehat{\varphi})),V (t-\eta(\widehat{\varphi}))) \over
f(T(t),V (t))}\right]
$$
$$ + e^{-\omega h} \left[ f(T(t-\eta(u_t)),V (t-\eta(u_t))) -  f(T(t-\eta(\widehat{\varphi})),V
(t-\eta(\widehat{\varphi})))\right] .
$$
To save the space we omit  long computations where we intensively
used equations (\ref{sdd-vir-stationary-2}), for example,
${e^{-\omega h}\over \delta + p\widehat{Y}} = {\widehat{T^{*}}\over
f(\widehat{T},\widehat{V})}$. Next, we add $\pm \left( 1 - {
V(t)\over \widehat{V} }\cdot {f({T}(t),\widehat{V})\over
f({T}(t),{V}(t))}\right) $ into $[...]$ to get
$${d\over dt}U^{1}(t) =  \left( 1- {f(\widehat{T},\widehat{V})\over f({T}(t),\widehat{V})}\right)
 e^{-\omega h} d \left( \widehat{T} -  T(t)\right)
$$
$$ - \widehat{T^{*}}(\delta + p\widehat{Y}) \left[ {f(\widehat{T},\widehat{V})\over f({T}(t),\widehat{V})}
+  { T^{*}(t)\cdot \widehat{V}\over \widehat{T^{*}}\cdot  V(t)} + {
V(t)\over \widehat{V} }\cdot {f({T}(t),\widehat{V})\over
f({T}(t),{V}(t))} + {\widehat{T^{*}}\over T^{*}(t)} \cdot  {
f(T(t-\eta(u_t)),V (t-\eta(u_t))) \over f(\widehat{T},\widehat{V})}
 \right. $$
$$ \left.  
- 4  - \ln { f(T(t-\eta(\widehat{\varphi})),V
(t-\eta(\widehat{\varphi}))) \over f(T(t),V (t))} + \left\{{
V(t)\over \widehat{V} } - {f({T}(t),{V}(t))\over
f({T}(t),\widehat{V})} + 1 - { V(t)\over \widehat{V} }\cdot
{f({T}(t),\widehat{V})\over f({T}(t),{V}(t))} 
 \right\} \right]
$$
$$ + e^{-\omega h} \left[ f(T(t-\eta(u_t)),V (t-\eta(u_t))) -  f(T(t-\eta(\widehat{\varphi})),V
(t-\eta(\widehat{\varphi})))\right].
$$
For the sum $\{...\}$ above, the particular form of the function $f$
(see (\ref{sdd-vir-03})) and computations give
\begin{equation}\label{sdd-vir-08}
{ V(t)\over \widehat{V} } - {f({T}(t),{V}(t))\over
f({T}(t),\widehat{V})} + 1 - { V(t)\over \widehat{V} }\cdot
{f({T}(t),\widehat{V})\over f({T}(t),{V}(t))}  = {(V(t)-
\widehat{V})^2 \, k_2(1+k_1 T(t))\over \widehat{V} (1+k_1 T(t)+k_2
\widehat{V})(1+k_1 T(t)+k_2 V(t))}.
\end{equation}
Now we add $\pm  {\widehat{T^{*}}\over T^{*}(t)} \cdot  {
f(T(t-\eta(\widehat{\varphi})),V (t-\eta(\widehat{\varphi}))) \over
f(\widehat{T},\widehat{V})}$ into $[...]$ above and substitute
(\ref{sdd-vir-08}) to obtain
$${d\over dt}U^{1}(t) =  \left( 1- {f(\widehat{T},\widehat{V})\over f({T}(t),\widehat{V})}\right)
 e^{-\omega h} d \left( \widehat{T} -  T(t)\right)
$$
$$ - \widehat{T^{*}}(\delta + p\widehat{Y}) \left[ {f(\widehat{T},\widehat{V})\over f({T}(t),\widehat{V})}
+  { T^{*}(t)\cdot \widehat{V}\over \widehat{T^{*}}\cdot  V(t)} + {
V(t)\over \widehat{V} }\cdot {f({T}(t),\widehat{V})\over
f({T}(t),{V}(t))} + {\widehat{T^{*}}\over T^{*}(t)} 
\cdot  {
f(T(t-\eta(\widehat{\varphi})),V (t-\eta(\widehat{\varphi}))) \over
f(\widehat{T},\widehat{V})}
 \right. $$
$$ \left.  
- 4  - \ln { f(T(t-\eta(\widehat{\varphi})),V
(t-\eta(\widehat{\varphi}))) \over f(T(t),V (t))} + {(V(t)-
\widehat{V})^2 \, k_2(1+k_1 T(t))\over \widehat{V} (1+k_1 T(t)+k_2
\widehat{V})(1+k_1 T(t)+k_2 V(t))} \right]
$$
\begin{equation}\label{sdd-vir-09}
 + e^{-\omega h} \left( 1 -{\widehat{T^{*}}\over T^{*}(t)} \right) \left[ f(T(t-\eta(u_t)),V
(t-\eta(u_t))) -  f(T(t-\eta(\widehat{\varphi})),V
(t-\eta(\widehat{\varphi})))\right].
\end{equation}

The first four terms in $[...]$ above suggest to split the logarithm
as follows
$$ \ln { f(T(t-\eta(\widehat{\varphi})),V
(t-\eta(\widehat{\varphi}))) \over f(T(t),V (t))}
$$
\begin{equation}\label{sdd-vir-10}
= \ln {f(\widehat{T},\widehat{V})\over f({T}(t),\widehat{V})} + \ln
{ T^{*}(t)\cdot \widehat{V}\over \widehat{T^{*}}\cdot  V(t)} + 
\ln \left( { V(t)\over \widehat{V} }\cdot
{f({T}(t),\widehat{V})\over f({T}(t),{V}(t))} \right) + 
\ln \left( {\widehat{T^{*}}\over T^{*}(t)} \cdot  {
f(T(t-\eta(\widehat{\varphi})),V (t-\eta(\widehat{\varphi}))) \over
f(\widehat{T},\widehat{V})} \right).
\end{equation}

Substitution of (\ref{sdd-vir-10}) into  (\ref{sdd-vir-09})
  implies
$${d\over dt}U^{1}(t) = \left( 1- {f(\widehat{T},\widehat{V})\over f({T}(t),\widehat{V})}\right)
 e^{-\omega h} d \left( \widehat{T} -  T(t)\right)  -  {(V(t)-
\widehat{V})^2 \, \cdot \widehat{T^{*}}(\delta +
 p\widehat{Y})\, k_2(1+k_1 T(t))\over \widehat{V} (1+k_1 T(t)+k_2
\widehat{V})(1+k_1 T(t)+k_2 V(t))}
$$
$$ - \widehat{T^{*}}(\delta + p\widehat{Y}) \left[
v \left( {f(\widehat{T},\widehat{V})\over
f({T}(t),\widehat{V})}\right) + v
\left( { T^{*}(t)\cdot \widehat{V}\over \widehat{T^{*}}\cdot  V(t)}\right) + 
v \left( { V(t)\over \widehat{V} }\cdot {f({T}(t),\widehat{V})\over
f({T}(t),{V}(t))} \right) \right.
$$
$$ \left. + \, 
v \left( {\widehat{T^{*}}\over T^{*}(t)} \cdot  {
f(T(t-\eta(\widehat{\varphi})),V (t-\eta(\widehat{\varphi}))) \over
f(\widehat{T},\widehat{V})} \right) \right]
$$
\begin{equation}\label{sdd-vir-11}
 + e^{-\omega h} \left( 1 -{\widehat{T^{*}}\over T^{*}(t)} \right) \left[ f(T(t-\eta(u_t)),V
(t-\eta(u_t))) -  f(T(t-\eta(\widehat{\varphi})),V
(t-\eta(\widehat{\varphi})))\right].
\end{equation}
Here we used the function $v(x)=x-1-\ln x$ to save the space. Next,
we can rewrite the first term in (\ref{sdd-vir-11}), using
(\ref{sdd-vir-03}),
$$ \left( 1- {f(\widehat{T},\widehat{V})\over f({T}(t),\widehat{V})}\right)
 e^{-\omega h} d \left( \widehat{T} -  T(t)\right) = - \left( T(t) -  \widehat{T}\right)^2
 { e^{-\omega h} d (1+k_2\widehat{V})\over {T}(t) (1 + k_1\widehat{T}+ k_2\widehat{V})}.
$$
We substitute the last equality into (\ref{sdd-vir-11}) to get
\begin{equation}\label{sdd-vir-12}
{d\over dt}U^{1}(t) = - D^{1}(t) + S^{1}(t),
\end{equation}
where 
$$ D^{1}(t)\equiv  \left( T(t) -  \widehat{T}\right)^2 \cdot
 { e^{-\omega h} d (1+k_2\widehat{V})\over {T}(t) (1 + k_1\widehat{T}+
 k_2\widehat{V})}
 +  {(V(t)-
\widehat{V})^2 \, \cdot \widehat{T^{*}}(\delta +
 p\widehat{Y})\, k_2(1+k_1 T(t))\over \widehat{V} (1+k_1 T(t)+k_2
\widehat{V})(1+k_1 T(t)+k_2 V(t))}
$$
$$ + \widehat{T^{*}}(\delta + p\widehat{Y}) \left[
v \left( {f(\widehat{T},\widehat{V})\over
f({T}(t),\widehat{V})}\right) + v
\left( { T^{*}(t)\cdot \widehat{V}\over \widehat{T^{*}}\cdot  V(t)}\right) + 
v \left( { V(t)\over \widehat{V} }\cdot {f({T}(t),\widehat{V})\over
f({T}(t),{V}(t))} \right) \right.
$$
\begin{equation}\label{sdd-vir-13}
 \left. + \, 
v \left( {\widehat{T^{*}}\over T^{*}(t)} \cdot  {
f(T(t-\eta(\widehat{\varphi})),V (t-\eta(\widehat{\varphi}))) \over
f(\widehat{T},\widehat{V})} \right) \right],
\end{equation}
\begin{equation}\label{sdd-vir-14}
 S^{1}(t) \equiv  e^{-\omega h} \left( 1 -{\widehat{T^{*}}\over T^{*}(t)} \right) \left[ f(T(t-\eta(u_t)),V
(t-\eta(u_t))) -  f(T(t-\eta(\widehat{\varphi})),V
(t-\eta(\widehat{\varphi})))\right].
\end{equation}
One can see, using $v(x)\ge 0$, that $D^{1}(t)\ge 0$.

{\bf Remark 9.} {\it It is easy to check  that $D^{1}(t)= 0$ if and
only if $T(t)=\widehat{T}, V(t)=\widehat{V},
T^{*}(t)=\widehat{T^{*}},$ $ f(T(t-\eta(\widehat{\varphi})),V
(t-\eta(\widehat{\varphi}))) = f(\widehat{T},\widehat{V})$. It
follows from the property $v(x)=0$ if and only if $x=1$.
}
\smallskip

 The sign of $S^{1}(t)$ is undefined. Our goal is to show that
$(-D^{1}(t)+ S^{1}(t))\le 0$, i.e. ${d\over dt}U^{1}(t)\le 0$ (see
(\ref{sdd-vir-12})) and ${d\over dt}U^{1}(t)=0$ at the stationary point only.
 To estimate the term $S^{1}(t)$ we notice that the functional
response $f(T,V)$, given by (\ref{sdd-vir-03}), is Lipschitz
$$ |f(T,V) - f(\widetilde{T},\widetilde{V})| \le L^f_1\cdot |T-\widetilde{T}|
+ L^f_2\cdot |V-\widetilde{V}|.
$$
It implies
$$\left| f(T(t-\eta(u_t)),V
(t-\eta(u_t))) -  f(T(t-\eta(\widehat{\varphi})),V
(t-\eta(\widehat{\varphi})))\right|
$$
$$\le L^f_1\cdot \left| T(t-\eta(u_t))-T(t-\eta(\widehat{\varphi}))\right|
+ L^f_2\cdot \left| V(t-\eta(u_t))- V
(t-\eta(\widehat{\varphi})))\right|.
$$

{\bf Remark 10.} {\it Both coordinates $T(t)$ and $V(t)$ of a solution
$u(t)$ of (\ref{sdd-vir-04}) are Lipschitz in time.  We denote the
corresponding Lipschitz constants for arbitrary solution as $L^T_u,
L^V_u$. It is easy to see that for any $\delta$-neighborhood of the
stationary point $\widehat \varphi$ the Lipschitz constants   of
arbitrary solution $|u_t - \widehat \varphi| \le \delta $ (inside of
the neighborhood) is uniformly bounded i.e. $L^T_u\le L^{T,\delta},
L^V\le L^{V,\delta}$. Moreover, $L^{T,\delta} \to 0, L^{V,\delta}\to
0$ as $\delta\to 0$.
}
\smallskip

We continue, using the assumptions on the state-dependent delay
$\eta$ (see (\ref{sdd-vir-05}) and (\ref{sdd-vir-06})),
$$ \left| f(T(t-\eta(u_t)),V (t-\eta(u_t))) -
f(T(t-\eta(\widehat{\varphi})),V (t-\eta(\widehat{\varphi})))\right|
$$
$$\le \left( L^f_1 L^{T,\delta} + L^f_2 L^{V,\delta}
\right) \cdot \left| \eta(u_t)-\eta(\widehat{\varphi})\right| \le
\left( L^f_1 L^{T,\delta} + L^f_2 L^{V,\delta} \right) \cdot c_\eta
\left( (T(t)-\widehat{T})^2 + (V(t)-\widehat{V})^2\right).
$$
This and (\ref{sdd-vir-14}) give the estimate
\begin{equation}\label{sdd-vir-15}
\left| S^{1}(t) \right| \le e^{-\omega h} \left| 1 -{\widehat{T^{*}}\over
T^{*}(t)} \right| \cdot  \left( L^f_1 L^{T,\delta} + L^f_2
L^{V,\delta} \right)\cdot c_\eta \left( (T(t)-\widehat{T})^2 +
(V(t)-\widehat{V})^2\right).
 \end{equation}
 Now we can choose small
enough $\delta$ (see remark)  to make the coefficient $e^{-\omega h}
\left| 1 -{\widehat{T^{*}}\over T^{*}(t)} \right| \cdot  \left(
L^f_1 L^{T,\delta} + L^f_2 L^{V,\delta} \right)\cdot c_\eta$ in
(\ref{sdd-vir-15}) arbitrary small, which implies (see the form of
$D^{1}(t)$ (\ref{sdd-vir-13})) the desired property ${d\over
dt}U^{1}(t) = - D^{1}(t) + S^{1}(t)< 0$.

{\bf Remark 11.} {\it It is easy to see from the calculations above
that the small value of $|T^{*}(t)-\widehat{T^{*}}|$ alone gives
${d\over dt}U^{1}(t) < 0$. No need to ask the values of
$|T(t)-\widehat{T}|, |V(t)-\widehat{V}|$ to be small. On the other
hand, the small values of $|T(t)-\widehat{T}|, |V(t)-\widehat{V}|$
give ${d\over dt}U^{1}(t) < 0$ without the need of small
$|T^{*}(t)-\widehat{T^{*}}|$.
}

The previous remark suggests that, alternatively, the smaller value
of constant $c_\eta$ (see (\ref{sdd-vir-06}) and (\ref{sdd-vir-15}))
the bigger the set where ${d\over dt}U^{1}(t) < 0$ holds. The latter
implies that in case $c_\eta$ is sufficiently small, the stationary
solution is {\tt globally} stable. Here the LaSalle invariance
principle (see e.g. \cite[Theorem 5.17, p.80]{Smith-2011-book}) is
applied. To show that the maximal invariant subset of ${d\over
dt}U^{1}(t)=0$ is the stationary solution only, we use the form of
${d\over dt}U^{1}(t)$ (we have shown that ${d\over dt}U^{1}(t)=0$
iff $D^{1}(t)=0$). We notice that $T(t)=\widehat{T},
T^{*}(t)=\widehat{T^{*}},V(t)=\widehat{V}$ follows immediately from
(\ref{sdd-vir-13}). Hence the second equation of (\ref{sdd-vir-04})
implies $0=\dot T^{*}(t) = e^{-\omega h} f(\widehat{T},\widehat{V})
- \delta \widehat{T^{*}} - p Y(t)\widehat{T^{*}}$. Similarly, the
third equation of (\ref{sdd-vir-04}) gives   $0=\dot V(t) = N\delta
\widehat{T^{*}} - c \widehat{V} - qA(t)\widehat{V}$. Clearly,
(\ref{sdd-vir-stationary-1}) implies $Y(t)\equiv \widehat{Y},
A(t)\equiv \widehat{A}$.
The proof of Theorem 8 is complete. $\Box $

\subsection{General case}


Now we are interested in  continuously differentiable solutions,
given by Theorem 6. As we mentioned above, for any solution $u$,
satisfying $u_0 \in \Omega_C$, one has $u_t \in \Omega_{\cal F} $
for $t>h$.

\medskip

{\bf Theorem 12}. {\it
Let assumptions (H2) and (H3) be satisfied (see (\ref{sdd-vir-H2}), (\ref{sdd-vir-H3})).
Assume the state-dependent delay $\eta : C \to [0,h]$ is
continuously differentiable  in a $\mu$-neighborhood of the
stationary solution $\widehat{\varphi}$ 
and satisfies $(H1_\eta)$ (see (\ref{sdd-vir-H1})).

Then the stationary solution
$\widehat{\varphi}=(\widehat{T},\widehat{T^{*}},\widehat{V},\widehat{Y},\widehat{A})$
of (\ref{sdd-vir-04}) is locally asymptotically stable. }

\smallskip

{\it Proof of Theorem 12.} Let us introduce the following Lyapunov
functional  with state-dependent delay along a solution of
(\ref{sdd-vir-04})
$$U^{sdd}(t) \equiv  \left( T(t) - \widehat T -
\int^{T(t)}_{\widehat T} {f(\widehat{T},\widehat{V})\over
f(\theta,\widehat{V})} \, d\theta  \right)e^{-\omega h}
 + \widehat{T^{*}}\cdot v\left({T^{*}(t)\over \widehat{T^{*}}}\right) +
{\delta + p\widehat{Y}\over N\delta} \cdot v\left({V(t)\over
\widehat V}\right) $$
\begin{equation}\label{sdd-vir-Lyapuniv-functional-2}
+ {p\over \beta} \cdot v\left({Y(t)\over \widehat Y}\right) + 
{q\over Ng}\left( 1+ {p\widehat{Y}\over \delta}\right) \cdot
v\left({A(t)\over \widehat A}\right) + (\delta + p\widehat{Y})
\widehat{T^{*}} \int^t_{t-\eta(u_t)}
v\left({f(T(\theta),V(\theta))\over f(\widehat{T},\widehat{V})}
\right) \, d\theta.
\end{equation}
A particular case of the constant delay functional was considered in
\cite{Zhao-Xu_EJDE-2014} (see (\ref{sdd-vir-Lyapuniv-functional-1})
above).
The difference is in the state-dependence of the lower bound of the
last integral in (\ref{sdd-vir-Lyapuniv-functional-2}).

Let us compute the time derivative of the last integral along a
continuously differentiable solution
$${d\over dt} \left( \int^t_{t-\eta(u_t)}
v\left({f(T(\theta),V(\theta))\over f(\widehat{T},\widehat{V})}
\right) \, d\theta \right) $$
$$= v\left({f(T(t),V(t))\over
f(\widehat{T},\widehat{V})} \right) -
v\left({f(T(t-\eta(u_t)),V(t-\eta(u_t)))\over
f(\widehat{T},\widehat{V})} \right) \cdot \left( 1 - {d\over
dt}\eta(u_t) \right)
$$
Comparing with the computations of ${d\over dt}U^{1}(t)$, we see the
main difference in the appearance of the term
\begin{equation}\label{sdd-vir-07}
S^{sdd}(t)\equiv
 - v\left({f(T(t-\eta(u_t)),V(t-\eta(u_t)))\over
f(\widehat{T},\widehat{V})} \right)  \cdot {d\over dt}\eta(u_t).
\end{equation}

\medskip

{\bf Remark 13.} {\it We notice that for any $u\in C^1([-h,b);R^5)$ one
has for $t\in [0,b)$
$${d\over dt}\eta(u_t) = [(D \eta)(u_t)](\dot u_t),$$
where $[(D \eta)(u_t)](\cdot )$ is the Frechet derivative of $\eta$
at point $u_t$. Hence, (for a solution in $\mu$-neighborhood of the
stationary solution $\widehat{\varphi}$) the estimate $|{d\over
dt}\eta(u_t)| \le || (D \eta)(u_t)||_{L(C;R)} \cdot ||\dot u_t
||_{C} \le  \mu\, || (D
\eta)(u_t)||_{L(C;R)}$ guarantees the property 
\begin{equation}\label{sdd-vir-30}
\left| {d\over dt}\eta(u_t) \right| \le \alpha_\mu \mbox{ with }
\alpha_\mu \to 0
 \mbox{ as } \mu\to 0.
\end{equation}
due to the boundedness of  $|| (D \eta)(\psi)||_{L(C;R)}$ as $\mu\to
0$ (here $||\psi - \widehat{\varphi}||_C < \mu$).
}

\smallskip

The time derivative of $U^{sdd}(t)$ along a continuously
differentiable solution $u$ of
(\ref{sdd-vir-04}) 
is computed similar to ${d\over dt}U^{1}(t)$ in the previous section.
We use (\ref{sdd-vir-stationary-2}) to get
$${d\over dt}U^{sdd}(t) = - D^{sdd}(t) + S^{sdd}(t),
$$
where 
$$ D^{sdd}(t) 
\equiv  \left( T(t) -  \widehat{T}\right)^2 \cdot
 { e^{-\omega h} d (1+k_2\widehat{V})\over {T}(t) (1 + k_1\widehat{T}+
 k_2\widehat{V})}
 +  {(V(t)-
\widehat{V})^2 \, \cdot \widehat{T^{*}}(\delta +
 p\widehat{Y})\, k_2(1+k_1 T(t))\over \widehat{V} (1+k_1 T(t)+k_2
\widehat{V})(1+k_1 T(t)+k_2 V(t))}
$$
$$ + \widehat{T^{*}}(\delta + p\widehat{Y}) \left[
v \left( {f(\widehat{T},\widehat{V})\over
f({T}(t),\widehat{V})}\right) + v
\left( { T^{*}(t)\cdot \widehat{V}\over \widehat{T^{*}}\cdot  V(t)}\right) + 
v \left( { V(t)\over \widehat{V} }\cdot {f({T}(t),\widehat{V})\over
f({T}(t),{V}(t))} \right) \right.
$$
\begin{equation}\label{sdd-vir-16}
 \left. + \, 
v \left( {\widehat{T^{*}}\over T^{*}(t)} \cdot  {
f(T(t-\eta(u_t)),V (t-\eta(u_t))) \over
f(\widehat{T},\widehat{V})} \right) \right],
\end{equation}
and $S^{sdd}(t)$ 
is defined in (\ref{sdd-vir-07}).

First, we observe that $D^{sdd}(t)$ and $D^{1}(t)$ have similar forms, while $S^{sdd}(t)$ and $S^{1}(t)$ are
essentially different. Moreover, in general case, we can not use (\ref{sdd-vir-05}), (\ref{sdd-vir-06}).

Our goal is to prove that there is a neighborhood of $\widehat u\in
C$, where ${d\over dt}U^{sdd}(t)< 0$ (except the point $\widehat u
$). We notice that $D^{sdd}(t)\ge 0$, while the sign of $S^{sdd}(t)$
is undefined. We will show that there is a neighborhood of the
stationary point, where $|S^{sdd}(t)| < D^{sdd}(t)$.

Let us consider the following auxiliary functionals $D^{(5)}(x)$ and
$S^{(5)}(x)$, defined on $R^5$, where we simplify notations
$x=(x^{(1)},x^{(2)},x^{(3)},x^{(4)},x^{(5)})\in R^5$ for $x^{(1)}=T,
x^{(2)}= T^{*}, x^{(3)}=V, x^{(4)} = T(t-\eta), x^{(5)} = V(t-\eta)$
$$
D^{(5)}(x)\equiv \left( {f(\widehat{T},\widehat{V})\over f(x^{(1)},
\widehat{V})} -1 \right)^2 + \left( {x^{(2)}\cdot\widehat{V})\over
\widehat{T^{*}}\cdot x^{(3)}} -1 \right)^2 + \left( {x^{(3)}\cdot
f(x^{(1)},\widehat{V})\over \widehat{V}\cdot f(x^{(1)},x^{(3)})} -1
\right)^2 $$
\begin{equation}\label{sdd-vir-D5}
+ \left( {\widehat{T^{*}}\cdot f(x^{(4)},x^{(5)})\over x^{(2)}\cdot
f(\widehat{T},\widehat{V})} -1 \right)^2 + c^{(1)}\cdot \left(
x^{(1)} - \widehat{T} \right)^2 + c^{(2)}\cdot \left( x^{(3)} -
\widehat{V} \right)^2, \quad c^{(1)}, c^{(2)} >0.
\end{equation}
\begin{equation}\label{sdd-vir-S5}
S^{(5)}(x)\equiv \alpha\cdot v\left( { f(x^{(4)},x^{(5)})\over
f(\widehat{T},\widehat{V})}\right), \quad \alpha\ge 0.
\end{equation}
The reason to consider functions $D^{(5)}(x)$ and $S^{(5)}(x)$ comes
from the property (\ref{sdd-vir-v-functional-2}) of the function
$v$. One sees that $D^{(5)}(x)=0$ if and only if $x=
\widehat{u}\equiv (\widehat{T},\widehat{T^{*}},\widehat{V},\widehat
Y,\widehat A)$. Now we change the coordinates in $R^5$ to the
spherical ones
\begin{equation}\label{sdd-vir-sph-coord}
\left\{
\begin{array}{l}
  x^{(1)}= \widehat{T} + r\cos \xi_4\, \cos \xi_3\, \cos \xi_2\,\cos \xi_1,\\
   x^{(2)}= \widehat{T^{*}} + r\cos \xi_4\, \cos \xi_3\, \cos \xi_2\,\sin \xi_1, \\
   x^{(3)}= \widehat{V} + r\cos \xi_4\, \cos \xi_3\, \sin \xi_2, \\
   x^{(4)}= \widehat{Y} + r\cos \xi_4\, \sin \xi_3, \\
   x^{(5)}= \widehat{A} + r\sin \xi_4, \qquad r\ge 0,\, \xi_1\in
   [0,2\pi), \, \xi_i\in [-\pi/2, \pi/2],\, i=2,...,5.
\end{array}
\right.
\end{equation}
One can check that the form of $D^{(5)}(x)$ (see (\ref{sdd-vir-D5}))
gives the multiplier $r^2$ in front of the sum, i.e.
$D^{(5)}(x)=r^2\cdot \Phi (r,\xi_1,...,\xi_5)$, where $\Phi
(r,\xi_1,...,\xi_5)$ is continuous and $\Phi (r,\xi_1,...,\xi_5)\neq
0$ for $r\neq 0$. The last property is proved, for example, assuming
the opposite $\Phi (r^0,\xi^0_1,...,\xi^0_5)= 0$ for $r^0\neq 0$,
which contradicts (\ref{sdd-vir-v-functional-2}). Hence, the
classical extreme value theorem (the Bolzano–Weierstrass theorem)
shows that the continuous $\Phi$ on a closed neighborhood of
$\widehat{u}$  has a minimum $\Phi_{min}>0$. It gives $D^{(5)}(x)\ge
r^2\cdot \Phi_{min}$.

Now the similar arguments for $S^{(5)}(x)$ shows that $|S^{(5)}(x)|
\le \alpha_\mu \cdot r^2$ where the constant $\alpha_\delta \to 0$
as $\mu\to 0$ (see (\ref{sdd-vir-30})). Finally, we can choose a
small enough $\mu>0$ to satisfy $ \alpha_\mu < \Phi_{min}$ which
proves that  ${d\over dt}U^{sdd}(t)\le - c r^2\cdot (\Phi_{min} -
\alpha_\mu)< 0$.

The proof of Theorem 12 is complete. $\Box $

\medskip

{\bf Remark 14.} {\it One sees that $S^{(5)}(x)$ depends on $x^{(4)},x^{(5)}$ only (\ref{sdd-vir-S5}).
On the other hand, the variables $x^{(4)},x^{(5)}$ are used in $D^{(5)}(x)$ in one term
$\left( {\widehat{T^{*}}\cdot f(x^{(4)},x^{(5)})\over x^{(2)}\cdot
f(\widehat{T},\widehat{V})} -1 \right)^2$ only. It is important to mention that the term
in $D^{(5)}(x)$ is not enough to bound $|S^{(5)}(x)|$ i.e.
\begin{equation}\label{sdd-vir-18}
|S^{(5)}(x)| \equiv \left| \alpha\cdot v\left( { f(x^{(4)},x^{(5)})\over
f(\widehat{T},\widehat{V})}\right) \right|\nleqslant \left( {\widehat{T^{*}}\cdot f(x^{(4)},x^{(5)})\over x^{(2)}\cdot
f(\widehat{T},\widehat{V})} -1 \right)^2.
\end{equation}
The sum of all terms in (\ref{sdd-vir-D5}) is needed to bound
$|S^{(5)}(x)|$. To see it, one should compare the sets where each
functional vanishes. Denote the zero-sets as $Z_{S^{(5)}}$ and
$Z_{rhs}$ (for the right-hand side of (\ref{sdd-vir-18})). Then one
sees that $Z_{S^{(5)}} \nsubseteq Z_{rhs}$. Moreover, in any
neighborhood of the point
$(x^{(2)},x^{(4)},x^{(5)})=(\widehat{T^{*}},\widehat{T},\widehat{V})\in
R^3$ one can find points where the right-hand side of
(\ref{sdd-vir-18}) is zero, while the the left-hand side is
positive. Clearly, the coordinates of such points should satisfy
$f(x^{(4)},x^{(5)}) \neq f(\widehat{T},\widehat{V}) $,
$\widehat{T^{*}}\cdot f(x^{(4)},x^{(5)})= x^{(2)}\cdot
f(\widehat{T},\widehat{V})$.
}



\medskip

{\bf Acknowledgments.}  This work was supported in part by GA CR under project 16-06678S.

\end{document}